# Evaluation of the Efficiency and Comparison of Different Numerical Differentiation Methods on Three Case Studies

Hamidreza Moradi[1], Hamideh Hossei[2]

[1]Department of Mechanical Engineering and Engineering Science, The University of North Carolina at Charlotte, Charlotte, North Carolina, USA[1]

[2]Department of Infrastructure and Environmental System, The University of North Carolina at Charlotte, North Carolina, USA[2]

**Abstract**

Without question regarding its pivotal significance, the computation of function derivatives carries substantial weight within a multitude of engineering and applied mathematical fields. These encompass optimization, the development of nonlinear control systems, and the assessment of noisy time signals, among others. In this study, we have chosen three illustrative cases: the logistic model for population dynamics, temperature variation within buildings, and the determination of market equilibrium prices. The primary objective is to assess the effectiveness of various numerical differentiation techniques and conduct a comparative analysis of the outcomes for each of these case studies. To achieve this objective, we employed three distinct numerical differentiation techniques: The Forward, Backward, and Centered Finite-Difference methods, each executed in two different levels of precision, totaling six variations. Our findings clearly indicate that, for the initial case study, the methods characterized by lower computational costs (specifically, the Forward and Backward Finite-Difference methods) yield superior outcomes. In contrast, for the second case study, the Centered Finite-Difference method delivers better results. In the case of the third case study, our results reveal that none of the methods produce estimations that meet acceptable standards. It is noteworthy that the empirical equations for each case study have been validated against previous literature.

[1] Corresponding Author.
Reza91moradi@gmail.com
[2] shossei5@charlotte.edu

## 1. Introduction

The modern landscape of engineering and scientific research frequently confronts intricate relationships among variables that defy facile linear modeling (Besharat Ferdosi and Porbashiri 2022; Myers et al. 2012; Najmi and Hu 2023). In the effort to navigate this complexity, numerical models have emerged as invaluable tools for tackling such challenges. Through the method of differentiation, numerical models have proven their worth by unraveling the intricate patterns concealed within complex data, enabling engineers and scientists to make accurate predictions, optimize processes, and facilitate data-driven decision-making (Singh and Sharma 2023). Additionally, the utility of curve fitting extends to parameter estimation and model validation, enabling the refinement of theoretical constructs and the validation of underlying assumptions, thus bolstering the dependability of engineering solutions and mitigating the costly specter of errors (Abbey and Meloy 2017). The ability to align empirical data with theoretical models empowers researchers to refine their designs and validate their concepts, enhancing the reliability and robustness of their solutions (Arabahmadi et al. 2023; Rezaeian et al. 2019). In the domains of signal processing and control systems, numerical differentiation offers the means to distill vital information from noisy datasets, enabling the separation of signal from noise and facilitating effective signal analysis and system control (Ramos and Polyakov 2015).

Three compelling case studies serve as the foundation of this endeavor. The first of these case studies is rooted in the domain of population dynamics, which can be applied on the dynamic change of all types of living animals or even a group of dynamic objects (Mahdavi, 2022a, Mahdavi, 2022b). Past studies have frequently harnessed various iterations of the exponential model to forecast population changes in this context, bearing the potential to significantly impact vital sectors like agriculture and demography (Kefayat, and Thill, 2025; Yang, Wu, and McMILLAN 1989; Sameer, Yusuf, and Ukafor 2022). Furthermore, the introduction of the logistic model, which incorporates three distinct types of data, opens the door to more precise population projections, thus advancing the accuracy of forecasting (Tu 1996; Mahdavi, 2024).

Moving onto the second case study, our exploration shifts to the realm of thermal comfort within educational and residential buildings, environments where people spend a significant portion of their lives. While previous research has made substantial contributions to this domain, a predominant focus has centered on specific educational facilities, and an important distinction regarding thermal comfort experienced during periods of occupancy versus non-occupancy has often been overlooked (Heracleous and Michael 2019; Lu et al. 2016; Shriram, K, and Swaminathan 2018; Branco et al. 2015; Quang et al. 2014). In this case study, we aim to examine how numerical models can effectively predict changes in the internal temperature of buildings, thereby addressing an essential aspect of occupant comfort and building energy management (Kermani et al., 2023).

The third and final case study shifts the focus to the domain of financial economics, specifically the equilibrium price of financial markets. This economic problem unfolds in two stages: first, the objective is to determine the optimal portfolio for each agent, considering the current prices of basic securities, followed by an examination of the equilibrium price of each basic security to ensure that the total supply matches the total demand. Within the confines of this third case study, our pursuit centers on refining the accuracy of relationships that represent the reference data under consideration, thus addressing critical issues related to market dynamics and pricing mechanisms.

At the core of these diverse case studies lies the pivotal role of numerical differentiation and integration, critical tools for effective data processing, particularly in the context of numerical differentiation (Babaey, and Faragardi, 2025;). Existing numerical differentiation and integration formulae, such as the Divided Difference formula, Newton-Cotes formula, and Gaussian quadrature rules, while valuable, harbor a common limitation: their foundation on polynomial functions implies that they excel primarily when applied to data conforming to smooth mathematical functions (Mirnajafizadeh et al., 2024). In practice, empirical data frequently deviate from these smooth functions and may even assume a form reminiscent of factual-like-derivative functions. Consequently, there arises a pressing need for numerical differentiation and integration techniques that can adapt to this variability, ensuring accurate results in diverse scenarios (Qu 1996).

The estimation of derivatives from numerical data, a classical challenge in the realm of data analysis, encounters a recurring issue: small errors in measurement data can yield significant errors in estimated derivatives, rendering the process ill-posed and challenging. Notably, one illustrative example of this challenge is the estimation of derivatives from noisy time signals, a problem of paramount significance and notorious complexity that has drawn extensive attention across various domains of engineering and applied mathematics (Mboup, Join, and Fliess 2007). This notoriety extends to the vast applicability of numerical differentiation, with its involvement in diverse numerical methods and essential processes such as the maximization of objective functions, the training of neural networks in machine learning (Widrow and Lehr 1990; Rastgoo et al., 2022; Moshiri, Heidari & Ghafoorzadeh-Yazdi, 2023; Ahmadi et al., 2023; Babaey and Ravindran, 2025a; Babaey, & Ravindran, 2025b;), and the sampling of high-dimensional posteriors in Bayesian models (Betancourt 2019; Neal 2011; Kucukelbir et al. 2017).

This research endeavor embarks on a comprehensive exploration of three distinct case studies, encompassing the logistic model of population change, temperature dynamics within buildings, and the equilibrium price of financial markets. These case studies are instrumental in evaluating and comparing the efficiency of various numerical differentiation methods, aiming to shed light on their efficacy and applicability in addressing the intricacies of real-world data. In the pursuit of this objective, the study not only scrutinizes the accuracy and stability of these methods but also delves into their computational efficiency, crucial for modern engineering and scientific endeavors where time and resources are of essence. Ultimately, this research seeks to offer a deeper understanding of numerical differentiation and its practical implications, enhancing the precision and reliability of results across diverse applications and domains.

## 2. Methods

In this study, we employ various numerical differentiation methods to analyze and compare their applicability and accuracy in three distinct case studies: population dynamics, temperature change of a building, and market equilibrium price. To begin, we gather the relevant data and formulate mathematical models for each case study, ensuring that these models represent the underlying dynamics accurately. Subsequently, we implement three commonly used finite difference numerical differentiation techniques, namely, central difference methods, forward and backward difference methods, to compute the derivatives required for our analysis. We then proceed to evaluate and compare the performance of these methods by quantifying their accuracy, stability, and computational efficiency in the context of each case study, providing valuable insights into the appropriate choice of differentiation technique for different real-world scenarios. We engage in three case studies aimed at investigating and conceptualizing three different categories of transformation: population change, temperature alterations, and shifts in market equilibrium prices. Each case study is rooted in established theories and empirical evidence and furnishing a comprehensive grasp of these dynamic occurrences.

## 2.1 Population Change

Consider a function, denoted as p(t), that represents the population at a specific time point, t. Despite the population being inherently discrete, it is generally large enough to make a reasonable assumption that p(t) exhibits characteristics of a continuous function. Our primary goal is to compute both the population's growth rate input and death rate output. The Malthusian model can be utilized to make these predictions. It aligns reasonably well with census data up to approximately 1900. However, beyond this point, the model overestimates the population, rendering it inadequate. Another commonly employed model for analyzing population dynamics is the Logistic Model, which is defined as follows:

$$\frac{dp}{dt} = -Ap(p - p_1) \qquad p(0) = p_0 \tag{1}$$

$$p(t) = \frac{p_0 p_1}{p_0 + (p_1 - p_0)e^{-Ap_1 t}} \tag{2}$$

In the academic context, starting with the reference year 1900 (t = 0) and an initial value of 76.09 ($p_0$), the task at hand is to ascertain the parameters A and $p_1$ within the given equation. To achieve this, the equation is adjusted to align with the available data, leading to the derived constants: A = -4.382e-07, $p_1$ = -2.921e+04, and the constant $p_0$ = 76.09. Ultimately, we have:

$$p(t) = \frac{-2.2226e + 06}{76.09 + (-2.9286e + 04)e^{-(0.0128)t}} \tag{3}$$

## 2.2 Temperature Change

The objective is to conduct an academic comparison of numerical models used to depict the 24-hour temperature variations within a structure, considering external temperature, internal heat generation, and the operation of heating or cooling systems. A common method for modeling indoor temperatures involves employing Compartmental Analysis. In this approach, we represent the building as a unified compartment, denoting the indoor temperature as T(t) at time t. The alteration in temperature over time is governed by all the elements contributing to heat generation or dissipation. Three key factors influencing the temperature inside a building. The first factor pertains to the heat generated by occupants, lighting, and equipment within the building, denoted as H(t), leading to a temperature increase. The second factor involves heating or cooling provided by the HVAC system (furnace or air conditioner), represented as U(t), and is typically expressed in terms of energy per unit time. The third factor considers the influence of external temperature, M(t). Consequently, we derive the following differential equation and temperature function:

$$\frac{dT}{dt} = K[M(t) - T(t)] + H(t) + U(t) \tag{4}$$

$$T(t) = B_2 - B_1F_1(t) + Ce^{(-k_1t)} \tag{5}$$

Where

$$F_1(t) = \frac{\cos(\omega t) + \left(\frac{\omega}{k_1}\right)\sin(\omega t)}{1 + \left(\frac{\omega}{k_1}\right)^2} \tag{6}$$

The constant C is selected when time is set to zero (t = 0), the temperature, denoted as T, is equivalent to the initial temperature, $T_0$. Thus, we have:

$$C = T_0 - B_2 + B_1F_1(0)$$

### 2.3 Market Equilibrium Price

Market equilibrium is a state in which the amount of a commodity demanded equals the amount supplied at a specific price. This analysis employs a linear model to investigate how market prices change in response to shifting conditions. The linear model for market equilibrium assumes that the demand and supply functions are represented as $q_d = d_0 - d_1p$ and $q_s = -s_0 + s_1p$, respectively, with p denoting the market price of the product, $q_d$ indicating the corresponding quantity demanded, $q_s$ representing the associated quantity supplied, and $d_0$, $d_1$, $s_0$, and $s_1$ being positive constants. These functional forms ensure that the fundamental principles of a downward-sloping demand curve and an upward-sloping supply curve are upheld. It can be readily demonstrated that the equilibrium price is calculated as $\dot{p} = \frac{(d_0+s_0)}{(d_1+s_1)}$.

Economists typically assume that markets are in equilibrium and justify this assumption with the help of stability arguments. We can consider a model that considers the expectations of agents. Based on the assumption that market demand and supply functions over time t >= 0 are as follows:

$$q_d(t) = d_0 - d_1 p(t) + d_2 \acute{p}(t) \; and \; q_s(t) = -s_0 - s_1 p(t) - s_2 \acute{p}(t) \tag{7}$$

respectively, where p(t) is the market price of the product, $q_d(t)$ is the associated quantity demanded, qs(t) is the associated quantity supplied, and $d_0$, $d_1$, $d_2$, $s_0$, $s_1$, and $s_2$ are all positive constants. Taking the ODE form as a starting point, we can derive the exact solution as follows:

$$p(t\acute{)} = \frac{\lambda a - \lambda b p}{1 - \lambda c} \tag{8}$$

$$P(t) = De^{\frac{\lambda b t}{c\lambda - 1}} + \frac{a}{b} \qquad \{a = d_0 + s_0 \; b = d_1 + s_1 \; c = d_2 + s_2\} \tag{9}$$

In the case of constants, they should be expressed as D = 3.282e-08, L = -7.314, a = 36.07, b = -0.01, c = 0.018. To improve the accuracy of our modeling, we explore a range of tailored equations that collectively provide a comprehensive portrayal of equilibrium price. This multifaceted methodology guarantees that our analysis is both resilient and representative of actual market equilibrium price change.

### 2.3 Numerical Differentiation

Work with three different numerical methods namely: Forward, and Centered finite-difference methods as shown in the equations 1-9, pave the way to estimate the derivative of a function at a particular point when an analytical expression for the function is not readily available. The forward method comes in handy when we are working with some data where we have information about the points after the current point. On the contrary, in backward method we work with points before the current point. In the centered method we need information of both previous and next point of the current point which require more data but result in more precise results. These methods can be classified into different orders, namely first-order and second order, based on the accuracy of the derivative estimation. First-order numerical differentiation methods involve estimating the derivative of a function using information from one neighboring point on either side of the evaluation point. Second-order numerical differentiation methods, on the other hand, use information from two neighboring points on both sides of the evaluation point. This study explores the forward, backward, and central numerical differentiation methods in both first and second orders.

**Forward Finite Difference**

In the forward difference method, the derivative is approximated using a point ahead of the evaluation point. The derivation of both first order and second order are shown below:

| First Derivative | Error |
|---|---|
| $f'(x_i) = \frac{f(x_{i+1}) - f(x_i)}{h}$ | *O(h)* |
| $f'(x_i) = \frac{-f(x_{i+2}) + 4f(x_{i+1}) - 3f(x_i)}{2h}$ | *O($h^2$)* |

**Backward Finite Difference**

In the backward difference method, a point behind the evaluation point is used. The derivation of both first order and second order are illustrated below:

| First Derivative | Error |
|---|---|
| $f'(x_i) = \frac{f(x_i) - f(x_{i-1})}{h}$ | *O(h)* |
| $f'(x_i) = \frac{3f(x_i) + 4f(x_{i-1}) + f(x_{i-2})}{2h}$ | *O($h^2$)* |

**Centered Finite Difference**

The central difference method is a second-order method and is considered more accurate than first-order methods. It computes the derivative by taking the difference between the function values at points both ahead and behind the evaluation point, which helps reduce the impact of noise and improves accuracy.

First Derivative Error

$$f'(x_i) = \frac{f(x_{i+1}) - f(x_{i-1})}{2h} \qquad O(h^2)$$

$$f'(x_i) = \frac{-f(x_{i+2}) + 8f(x_{i+1}) - 8f(x_{i-1}) + f(x_{i-2})}{12h} \qquad O(h^4)$$

Higher order differentiation formulas can be calculated using Taylor series expansion as illustrated below:

$$f(x_{i+1}) = f(x_i) + f(\acute{x}_i)h + \frac{f''(x_i)}{2}h^2 + \cdots$$

$$f(\acute{x}_i) = \frac{f(x_{i+1}) - f(x_i)}{h} - \frac{f''}{2}h + O(h^2)$$

$$f''(x_i) = \frac{f(x_{i+2}) - 2f(x_{i+1}) + f(x_i)}{h^2} + O(h)$$

$$f(\acute{x}_i) = \frac{f(x_{i+1}) - f(x_i)}{h} - \frac{f(x_{i+2}) - 2f(x_{i+1}) + f(x_i)}{2h^2}h + O(h^2)$$

$$f(\acute{x}_i) = \frac{-2f(x_{i+2}) + 4f(x_{i+1}) - 3f(x_i)}{2h} + O(h^2)$$

For each case study error has been calculated for each numerical differentiation method versus both empirical relationship and experimental. Consider that the data resulted from empirical equation can be shown as $E_1$, $E_2$, $E_3$… $E_n$, and the data from the experimental work is shown by $P_1$, $P_2$, $P_3$… $P_n$ as well as the data form each numerical differentiation method as $X_1$, $X_2$, $X_3$… $X_n$ we can calculate the error for each method as

$$Error_{wrt_Exp} = \frac{\sum_{i=1}^{n}(P_i - X_i)}{\sum_{i=1}^{n}(P_i)}$$

$$Error_{wrt_Emp} = \frac{\sum_{i=1}^{n}(E_i - X_i)}{\sum_{i=1}^{n}(E_i)}$$

## 3. Results and Discussion

The choice between forward, backward, and central numerical differentiation methods depend on the specific application, the level of accuracy required, and the nature of the data or function being analyzed. Comparing these methods, first-order methods are simple to implement but can be less accurate, particularly when dealing with noisy data or functions with rapid changes. Second-order methods, like the central difference method, provide more accurate derivative estimates because they consider points on both sides of the evaluation point, but they are slightly more complex to apply.

**3.1 Population dynamics**

In the first case study we see that the results are so close to each other and the best method to select in this case is the one with the least solution cost and time. Because both forward and backward difference methods need less information than other methods and are of the same precision order, it is assumed that both are the best choices in this case study. The error presented in Table1 indicates the difference between the types of the solver and empirical/experimental data.

When assessing the numerical differentiation model from Table 1, it is evident that Forward2 method exhibited the highest error, with empirical error measuring 3.3252. In contrast, this method demonstrated slightly better accuracy against the experimental error measuring 0.0977. The Backward2 and Centered2 method performed with similar accuracy, yielding experimental errors below 0.05. However, the Backward2 solver stood out with the lowest error at 0.0208, signifying its superior accuracy in approximating the population dynamics using the exponential model. This analysis highlights the varying levels of precision among these solvers, with ode45 delivering the most reliable results in this specific context.

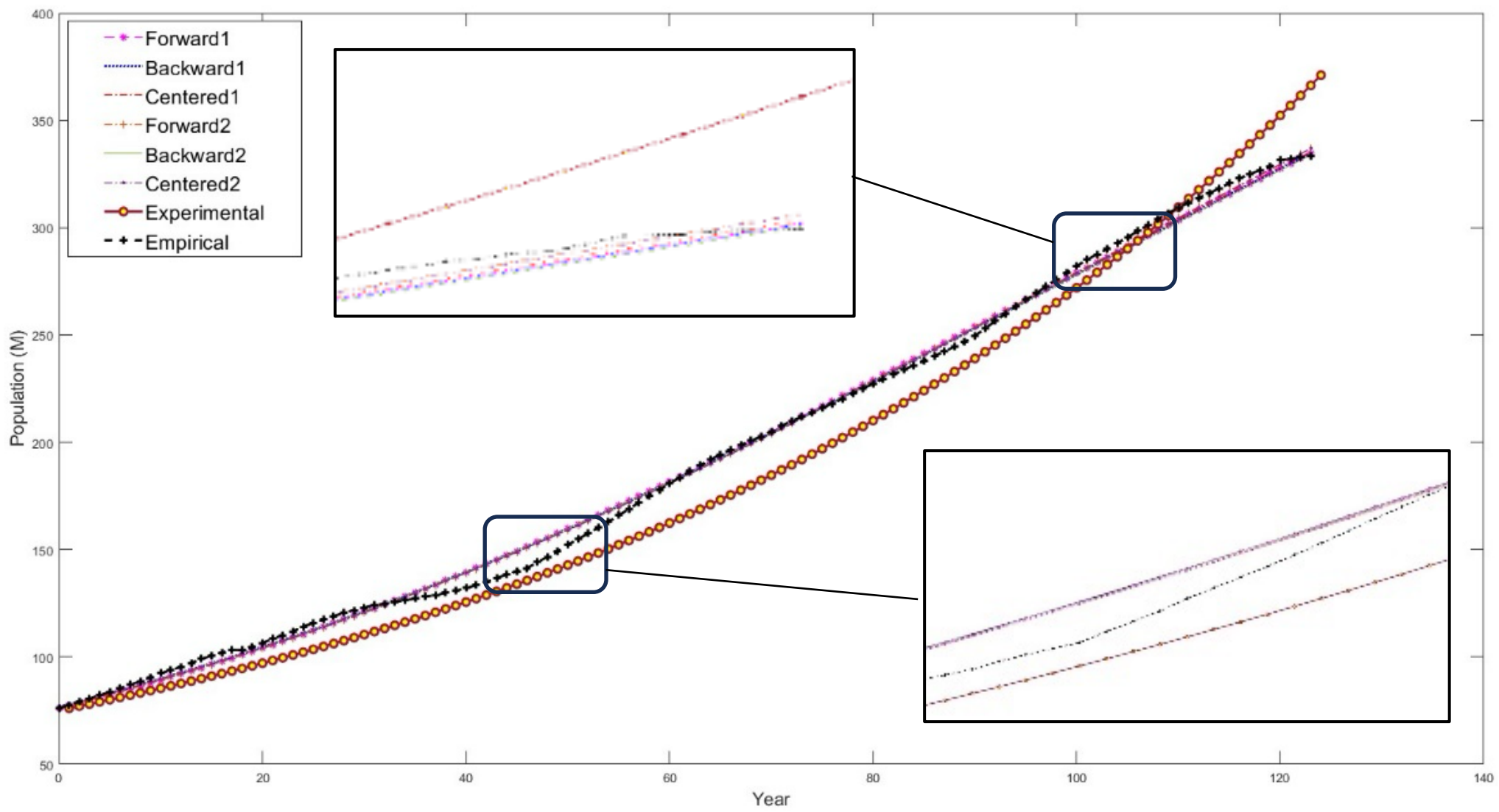

Figure 1. Comparison of different numerical differentiation methods for population dynamics

Table 1. Error values for population dynamics compared to both empirical and experimental results

| Error | Value |
|---|---|
| $e_{exp}^{forward1}$ | 0.0909 |
| $e_{emp}^{forward1}$ | 3.3182 |
| $e_{exp}^{backward1}$ | 0.0266 |
| $e_{emp}^{backward1}$ | 3.2519 |
| $e_{exp}^{centered1}$ | 0.0561 |
| $e_{emp}^{centered1}$ | 3.2823 |
| $e_{exp}^{forward2}$ | 0.0977 |
| $e_{emp}^{forward2}$ | 3.3252 |
| $e_{exp}^{backward2}$ | 0.0208 |
| $e_{emp}^{backward2}$ | 3.2460 |
| $e_{exp}^{centered2}$ | 0.0483 |
| $e_{emp}^{centered2}$ | 3.2743 |

### 3.2 Temperature Change

For the second case study we explicitly see the differences of each method. Therefore, we can choose the best method regarding the relative percentage errors associated with each numerical differentiation method. In this case study, the Forward differentiation error is higher at approximately 1.2108. In the empirical data, the error has also increased to 1.0046, but it remains lower than the error in the expected data, indicating that the empirical measurements are relatively more accurate for this method. For the Backward differentiation method, the error is even higher at 1.8456. The empirical data also has a higher error of 1.6407 but is more accurate than the expected data. For the Centered differentiation method, the error has increased to 1.2319, indicating decreased accuracy. In the empirical data, the error has also increased to 1.0258 but remains lower than the error in the expected data, suggesting that the empirical measurements are relatively more accurate for this method. In the second-order differentiation method, for the forward method, the error is relatively high at 1.2888, indicating less accuracy. The empirical data also has a higher error of 1.0828 but is more accurate than the expected data. For the backward method the error is now 1.1240 and indicates relatively lower accuracy. The empirical data has a lower error of 0.9176 and is more accurate. The centered method for second-order differentiation shows a significantly lower error with a value of 0.4377, suggesting a more accurate result and the empirical data also exhibits a lower error of 0.2299 and is more accurate than the expected data. As the results indicate, both the numerical differentiation methods and the empirical data exhibit higher error values compared to the previous set of data. However, the second-order Centered solver against empirical data stood out with the lowest error term.

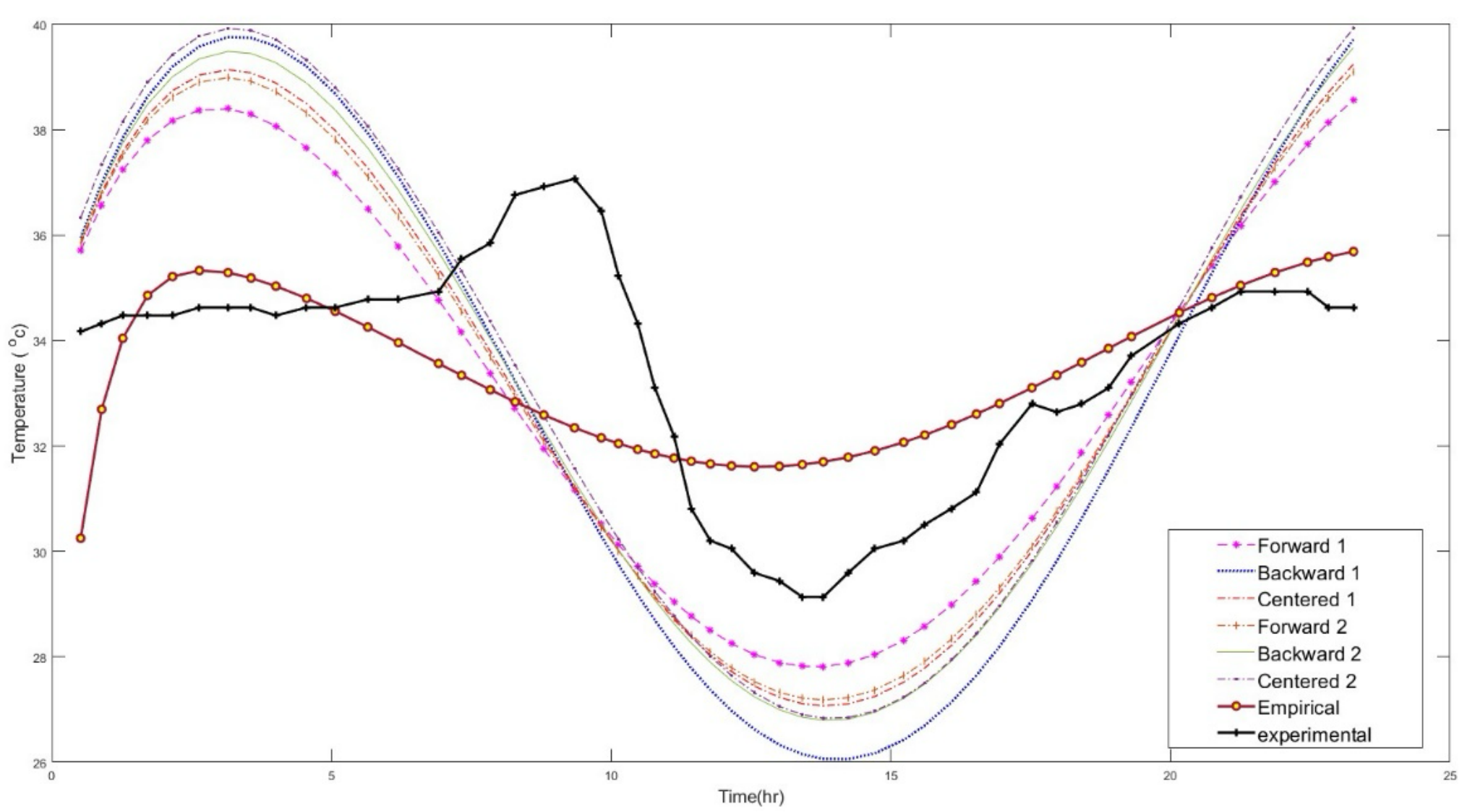

Figure 2. Comparison of different numerical differentiation methods for temperature change

Table 2. Error values for temperature change compared to both empirical and experimental results

| Error | Value |
| --- | --- |
| $e_{exp}^{forward1}$ | 1.2108 |
| $e_{emp}^{forward1}$ | 1.0046 |
| $e_{exp}^{backward1}$ | 1.8456 |
| $e_{emp}^{backward1}$ | 1.6407 |
| $e_{exp}^{centered1}$ | 1.2319 |
| $e_{emp}^{centered1}$ | 1.0258 |
| $e_{exp}^{forward2}$ | 1.2888 |
| $e_{emp}^{forward2}$ | 1.0828 |
| $e_{exp}^{backward2}$ | 1.1240 |
| $e_{emp}^{backward2}$ | 0.9176 |
| $e_{exp}^{centered2}$ | 0.4377 |
| $e_{emp}^{centered2}$ | 0.2299 |

### 3.3 Market Equilibrium Price

As we can see in figure 3 the nature of the market equilibrium price has a sharp slope after a certain amount of time. For this reason, each of the numerical differentiation methods end up with unusual rising behavior in the proximity of that area. This is an example of the cases where none of the numerical differentiation methods gives a reasonable answer to the equations. Table 3 presents a collection of error values associated with numerical differentiation techniques in comparison to empirical or experimental errors. The forward and backward difference methods yield relatively higher errors, as indicated by values around 34 to 121, while the centered difference method generally results in smaller errors, with values ranging from approximately 0.59 to 8.28. A general assertion resulted from table 3 is that for the phenomena with fluctuation behavior and cases with high slope of changes, numerical differentiation methods cannot converge to an acceptable answer. As a future work, the numerical differentiation methods with higher precision for which more information around the point of interest are needed can be used. Other methods to be used in these cases can be using ODE solvers directly rather than solving the equation with numerical differentiation methods.

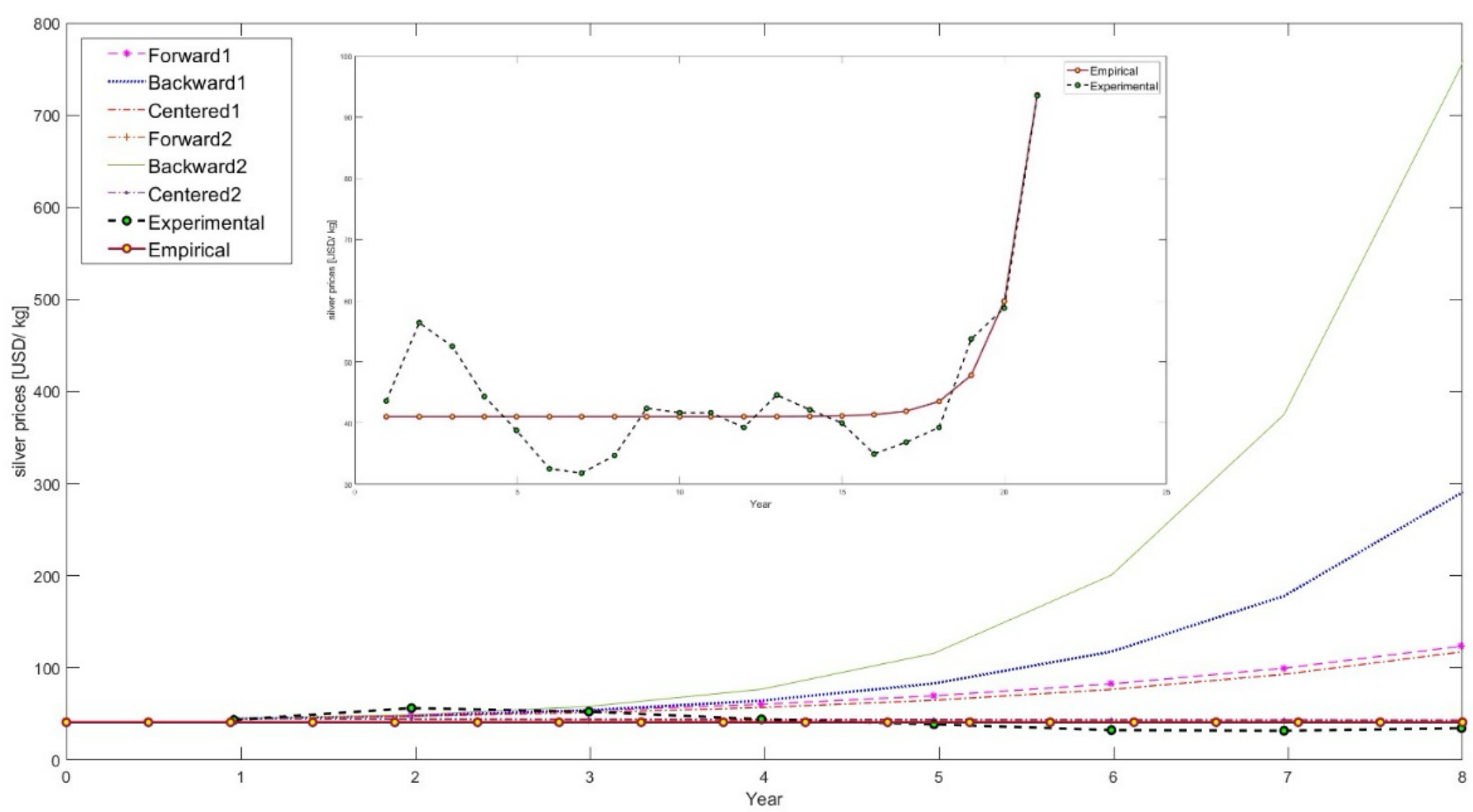

Figure 3. Comparison of different numerical differentiation methods for market equilibrium price

Table 3. Error values for market equilibrium price compared to both empirical and experimental results

| Error | Value |
|---|---|
| $e_{exp}^{forward1}$ | 34.2569 |
| $e_{emp}^{forward1}$ | 46.2297 |
| $e_{exp}^{backward1}$ | 53.7693 |
| $e_{emp}^{backward1}$ | 67.4821 |
| $e_{exp}^{centered1}$ | 27.7075 |
| $e_{emp}^{centered1}$ | 39.0962 |
| $e_{exp}^{forward2}$ | 2.7121 |
| $e_{emp}^{forward2}$ | 5.9638 |
| $e_{exp}^{backward2}$ | 103.4304 |
| $e_{emp}^{backward2}$ | 121.5720 |
| $e_{exp}^{centered2}$ | 0.5882 |
| $e_{emp}^{centered2}$ | 8.2771 |

## 4. Conclusion

In this work three case studies namely, logistic model of population change, temperature change in buildings and market equilibrium price have been considered to study the efficiency of the numerical differentiation methods and compare their results with each other. The results of the first case study show that there are small errors for each of the numerical differentiation methods so that can just consider the solution cost of each method for its superiority. Although the errors are small, the best method in this case is higher order backward finite-difference. For the second case study, we see there is an obvious difference between each method and based on the relative percentage error results the best numerical differentiation method in this case is higher order centered finite-difference. Due to the nature of the third case study for which we have a steep slope in the behavior of the results, all the numerical differentiation methods fail to estimate after a certain amount of time.